# A GEOMETRIC CHARACTERIZATION OF $C$-OPTIMAL DESIGNS FOR HETEROSCEDASTIC REGRESSION

By Holger Dette[1,2] and Tim Holland-Letz

*Ruhr University, Bochum*

We consider the common nonlinear regression model where the variance, as well as the mean, is a parametric function of the explanatory variables. The $c$-optimal design problem is investigated in the case when the parameters of both the mean and the variance function are of interest. A geometric characterization of $c$-optimal designs in this context is presented, which generalizes the classical result of Elfving [*Ann. Math. Statist.* **23** (1952) 255–262] for $c$-optimal designs. As in Elfving's famous characterization, $c$-optimal designs can be described as representations of boundary points of a convex set. However, in the case where there appear parameters of interest in the variance, the structure of the Elfving set is different. Roughly speaking, the Elfving set corresponding to a heteroscedastic regression model is the convex hull of a set of ellipsoids induced by the underlying model and indexed by the design space. The $c$-optimal designs are characterized as representations of the points where the line in direction of the vector $c$ intersects the boundary of the new Elfving set. The theory is illustrated in several examples including pharmacokinetic models with random effects.

**1. Introduction.** Nonlinear regression models are widely used to describe the relation between several variables [see Seber and Wild (1989), Ratkowsky (1983, 1990)]. Because of the broad applicability of such models in many fields, the problem of constructing optimum experimental designs for these models has found considerable interest in the literature. Early work has been

Received November 2008; revised February 2009.

[1]Supported by the Deutsche Forschungsgemeinschaft: Sonderforschungsbereich 475, Komplexitätsreduktion multivariater Datenstrukturen (Teilprojekt A2); Sonderforschungsbereich 823, Statistik nichlinearer dynamischer Prozesse.

[2]Supported in part by a NIH Grant IR01GM072876 and the BMBF-Grant SKAVOE.

*AMS 2000 subject classification.* 62K05.

*Key words and phrases.* $c$-optimal design, heteroscedastic regression, Elfving's theorem, pharmacokinetic models, random effects, locally optimal design, geometric characterization.

done by Chernoff (1953) who introduced the concept of locally optimal designs. These designs require an initial guess of the unknown parameters of the model and are used as benchmarks for many commonly used designs. Locally optimum designs for nonlinear models have been discussed by numerous authors [see Ford, Torsney and Wu (1992), Box and Lucas (1959), Haines (1993), Haines (1995), Biedermann, Dette and Zhu (2006), López-Fidalgo and Wong (2002), Dette, Melas and Pepelyshev (2004) among many others]. Most of the relevant literature discusses the design problem under the additional assumption of a constant variance, but much less work has been done for models with an heteroscedastic error structure. Wong and Cook (1993) studied $G$-optimal designs for models, when heteroscedasticity is present in the data. A systematic approach to optimal design problems for heteroscedastic linear models was given by Atkinson and Cook (1995), who derived the necessary information matrices in the case where the variance, as well as the mean, depend on the parameters of the model and the explanatory variables. For other work on optimal designs for specific regression models with heteroscedastic errors and various optimality criteria, we refer to King and Wong (1998), Ortiz and Rodríguez (1998), Fang and Wiens (2000), Brown and Wong (2000), Montepiedra and Wong (2001) and Atkinson (2008) among others.

The present paper is devoted to the local $c$-optimality criterion which determines the design such that a linear combination of the unknown parameters (specified by the vector $c$) has minimal variance [see, e.g., Pukelsheim (1981), Pázman and Pronzato (2009)]. Under the assumption of homoscedasticity, there exists a beautiful geometric characterization of locally $c$-optimal designs which is due to Elfving (1952) who considered the case of a linear model. Roughly speaking, the characterization of Elfving (1952) is possible because the Fisher information can be represented in the form

$$I(x,\theta) = f(x,\theta)f^T(x,\theta), \tag{1.1}$$

where $f(x,\theta)$ denotes a vector of functions corresponding to the particular model under consideration and $\theta$ is the vector of unknown parameters. The $c$-optimal design with support points $x_1,\dots,x_m$ and weights $p_1,\dots,p_m$ can be characterized as a representation of the point

$$\lambda c = \sum_{r=1}^{m} p_r \varepsilon_r f(x_r,\theta), \tag{1.2}$$

where $\varepsilon_1,\dots,\varepsilon_k \in \{-1,1\}$ and $\lambda > 0$ denotes a scaling factor such that the point $\lambda c$ is a boundary point of the Elfving set

$$\mathcal{R}_1 = \operatorname{conv}\{\varepsilon f(x,\theta) \mid x \in \mathcal{X}, \varepsilon \in \{-1,1\}\}. \tag{1.3}$$

Here, $\mathcal{X}$ denotes the design space, $\operatorname{conv}(\mathcal{A})$ is the convex hull of a set $\mathcal{A} \subset \mathbb{R}^p$. This result has been applied by numerous authors to derive $c$-optimal designs in linear and nonlinear regression models [see, e.g., Studden (1968), Han and Chaloner (2003), Ford, Torsney and Wu (1992), Chernoff and Haitovsky (1990), Fan and Chaloner (2003), Dette et al. (2008) among many others]. For a review on Elfving's theorem, we also refer to Fellmann (1999) and to the recent work of Studden (2005). On the other hand, in the case of heteroscedasticity, where both the mean and the variance depend on the explanatory variables and parameters of interest, the information matrix is usually of the form

$$I(x,\theta) = \sum_{j=1}^{k} f_j(x,\theta) f_j^T(x,\theta) \tag{1.4}$$

[see Atkinson and Cook (1995)] with $k \geq 2$. Here, $f_1, \dots, f_k$ represent certain features of the nonlinear regression model under consideration (see Section 2 for more details) and a geometric characterization of the Elfving type is not available. It is the purpose of the present paper to fill this gap and to provide a useful characterization of $c$-optimal designs for regression models with Fisher information of the form (1.4).

In Section 2, we introduce the basic notation of a heteroscedastic nonlinear regression model and demonstrate that the Fisher information matrix has, in fact, the form (1.4). Section 3 contains our main results. We derive a general equivalence theorem for locally $c$-optimal designs in nonlinear regression models with Fisher information matrix given by (1.4). These results are used to derive a geometric characterization of locally $c$-optimal designs which generalize the classical result of Elfving to regression models, where the mean as well as the variance depend on the explanatory variables and the parameters of interest. The corresponding Elfving space has the same dimension as a classical Elfving space but consists of the convex hull of a family of ellipsoids induced by the underlying model and indexed by the design space (see Section 3 for details). In the special case $k = 1$, the generalized Elfving set reduces to the classical set considered by Elfving (1952) and a similar result is obtained if the variance of the regression is a function of the mean. Finally, in Section 4, we illustrate the geometric characterization in several examples. In particular, we demonstrate the applicability of the results in the Michaelis–Menten model with exponentially decreasing variance and in a nonlinear random effect model used in toxicokinetics, where the Fisher information has a similar structure as in (1.4).

## 2. $c$-optimal designs for heteroscedastic regression models.

Consider the common nonlinear regression model where at a point $x$ a response $Y$ is observed with expectation

$$\mathrm{E}[Y|x] = \mu(x,\theta) \tag{2.1}$$

and variance

$$\operatorname{Var}[Y|x] = \sigma^2(x,\theta). \tag{2.2}$$

Here, $\theta \in \mathbb{R}^p$ denotes the vector of unknown parameters. Note that we do not exclude the case, where $\mu$ and $\sigma^2$ depend on different subsets of the parameter vector, that is $\theta = (\theta_1, \theta_2)$, $\mu(x,\theta) = \mu(x,\theta_1)$; $\sigma^2(x,\theta) = \sigma^2(x,\theta_2)$. We assume that $n$ independent observations $Y_1, \ldots, Y_n$ are available under experimental conditions $x_1, \ldots, x_n \in \mathcal{X}$, where $\mathcal{X}$ denotes the design space. We define $\mu_\theta = (\mu(x_1,\theta), \ldots, \mu(x_n,\theta))^T$ as the vector of the expected responses, $\sigma_\theta^2 = (\sigma^2(x_1,\theta), \ldots, \sigma^2(x_n,\theta))^T$ as the vector of variances and

$$\Sigma_\theta = \operatorname{diag}(\sigma^2(x_1,\theta), \ldots, \sigma^2(x_n,\theta)) \tag{2.3}$$

as the covariance matrix of the random vector $Y = (Y_1, \ldots, Y_n)^T$. Under the additional assumption of a normal distribution, the Fisher information of $Y$ is given by the $p \times p$ matrix

$$I = \frac{d\mu_\theta^T}{d\theta} \Sigma_\theta^{-1} \frac{d\mu_\theta}{d\theta} + \frac{1}{2} \frac{d\sigma_\theta^2}{d\theta}^T \Sigma_\theta^{-2} \frac{d\sigma_\theta^2}{d\theta}. \tag{2.4}$$

An approximate experimental design is a discrete probability measure with masses $w_1, \ldots, w_m$ at points $x_1, \ldots, x_m \in \mathcal{X}$. These points define the distinct experimental conditions at which observations have to be taken and $w_1, \ldots, w_m > 0$, $\sum_{j=1}^m w_j = 1$ are positive weights representing the proportions of total observations taken at the corresponding points [see Silvey (1980), Atkinson and Donev (1992), Pukelsheim (1993), Randall, Donev and Atkinson (2007)]. If $N$ observations can be taken, a rounding procedure is applied to obtain integers $r_j$ from the not necessarily integer valued quantities $w_j N$ $(j = 1, \ldots, m)$ [see Pukelsheim and Rieder (1992)]. If the assumption of a normal distribution is made the analogue of the Fisher information matrix for an approximate design is the matrix

$$M(\xi,\theta) = \int_{\mathcal{X}} I(x,\theta)\, d\xi(x) \in \mathbb{R}^{p\times p}, \tag{2.5}$$

where

$$\begin{aligned} I(x,\theta) = {} & \frac{1}{\sigma^2(x,\theta)} \left(\frac{\partial \mu(x,\theta)}{\partial\theta}\right)^T \frac{\partial\mu(x,\theta)}{\partial\theta} \\ & + \frac{1}{2\sigma^4(x,\theta)} \left(\frac{\partial\sigma^2(x,\theta)}{\partial\theta}\right)^T \frac{\partial\sigma^2(x,\theta)}{\partial\theta} \end{aligned} \tag{2.6}$$

denotes the Fisher information at the point $x$ [see also Atkinson and Cook (1995) or Atkinson (2008)].

Under some additional assumptions of regularity, it can be shown that the asymptotic covariance matrix of the maximum likelihood estimate for the parameter $\theta$ is given by

$$\frac{1}{N} M^{-1}(\xi, \theta) \tag{2.7}$$

[see Jennrich (1969)]. An optimal design for estimating the parameter $\theta$ minimizes an appropriate function of this matrix, and numerous criteria have been proposed for this purpose [see Silvey (1980), Pukelsheim (1993) or Randall, Donev and Atkinson (2007)]. In this paper, we consider the local $c$-optimality criterion, which determines the design $\xi$ such that the expression

$$c^T M^{-}(\xi, \theta) c, \tag{2.8}$$

is minimal, where $c \in \mathbb{R}^p$ is a given vector, the minimum is calculated among all designs for which $c^T\theta$ is estimable, that is $c \in \text{Range}(M(\xi, \theta))$, and $M^{-}(\xi, \theta)$ denotes a generalized inverse of the matrix $M(\xi, \theta)$. The expression (2.8) is approximately proportional to the variance of the maximum likelihood estimate for the linear combination $c^T\theta$. It is shown in Pukelsheim (1993) that the expression (2.8) does not depend on the specific choice of the generalized inverse if the vector $c$ is estimable by the design $\xi$. Note that the criterion (2.8) is a local optimality criterion in the sense that it requires the specification of the unknown parameter $\theta$ [see Chernoff (1953)]. Locally optimal designs are commonly used as benchmarks for given designs used in applications. Moreover, in some cases, knowledge about the parameter $\theta$ is available from previous studies and locally $c$-optimal designs are robust with respect to misspecifications of $\theta$ [see Dette et al. (2008)]. For robustifications of locally optimal designs, the interested reader is referred to the work of Chaloner and Verdinelli (1995), Dette (1995) or Müller and Pázman (1998) among many others.

Note that the Fisher information (2.4) is of the form

$$I(x, \theta) = \sum_{\ell=1}^{k} f_\ell(x, \theta) f_\ell^T(x, \theta), \tag{2.9}$$

where $k = 2$ and

$$\begin{aligned} f_1(x, \theta) &= \frac{1}{\sigma(x, \theta)} \left( \frac{\partial \mu(x, \theta)}{\partial \theta} \right)^T, \\ f_2(x, \theta) &= \frac{1}{\sqrt{2}\sigma^2(x, \theta)} \left( \frac{\partial \sigma^2(x, \theta)}{\partial \theta} \right)^T. \end{aligned} \tag{2.10}$$

Because the formulation of our main results does not yield any additional complication in the subsequent discussion, we consider in the following models with a Fisher information of the form (2.9). We call a design $\xi_c$ minimizing $c^T M^-(\xi,\theta)c$ a locally $c$-optimal for a model with Fisher information matrix of the form (2.9), where the minimum is taken over the class of all designs for which $c^T\theta$ is estimable. The corresponding statements for the heteroscedastic nonlinear regression models are then derived as special cases.

**3. Elfving's theorem for heteroscedastic models.** The following result allows an easy verification of $c$-optimality for a given design.

THEOREM 3.1. *A design $\xi_c$ is locally $c$-optimal in a model with Fisher information matrix of the form (2.9) if and only if there exists a generalized inverse $G$ of the matrix $M(\xi_c,\theta)$, such that the inequality*

$$\sum_{\ell=1}^{k} \frac{(c^T G f_\ell(x))^2}{c^T M^-(\xi_c,\theta)c} \le 1 \tag{3.1}$$

*holds for all $x \in \mathcal{X}$. Moreover, there is equality in (3.1) at any support point of the design $\xi_c$.*

PROOF. Let $\Xi$ denote the set of all approximative designs on $\mathcal{X}$ and for fixed $\theta$ let

$$\mathcal{M} = \{M(\xi,\theta) \mid \xi \in \Xi\} \subset \mathbb{R}^{p\times p}$$

denote the set of all information matrices of the form (2.5), where $I(x,\theta) = \sum_{\ell=1}^k f_\ell(x,\theta) f_\ell^T(x,\theta)$. $\mathcal{M}$ is obviously convex and the information matrix $M(\xi_c,\theta)$ of a locally $c$-optimal design for which the linear combination $c^T\theta$ is estimable [i.e., $c \in \mathrm{Range}(M(\xi,\theta))$] maximizes the function $(c^T M^- c)^{-1}$ in the set $\mathcal{M} \cap \mathcal{A}_c$, where $\mathcal{A}_c = \{M(\xi,\theta)) \in \mathcal{M} | c \in \mathrm{Range}(M(\xi,\theta))\}$. Consequently, it follows from Theorem 7.19 in Pukelsheim (1993) that the design $\xi_c$ is $c$-optimal if and only if there exists a generalized inverse, say $G$, of the matrix $M(\xi_c,\theta)$ such that the inequality

$$\mathrm{tr}(AGcc^T G^T) \le c^T M^-(\xi_c,\theta)c$$

holds for all $A \in \mathcal{M}$, where there is equality for any matrix $A \in \mathcal{M}$ which maximizes $(c^T M^- c)^{-1}$ in the set $\mathcal{M}$. Note that the family $\mathcal{M}$ is the convex hull of the set

$$\left\{ \sum_{\ell=1}^{k} f_\ell(x,\theta) f_\ell^T(x,\theta) \,\middle|\, x \in \mathcal{X} \right\},$$

and therefore the assertion of Theorem 3.1 follows by a standard argument of optimal design theory [see, e.g., Silvey (1980)]. □

The following result gives the corresponding statement for the case of the nonlinear heteroscedastic regression model introduced in Section 2. It is a direct consequence of Theorem 3.1 observing that in this case we have $k = 2$ and the two functions $f_1$ and $f_2$ are given by (2.10).

COROLLARY 3.2. *Under the assumption of a normally distributed error, a design $\xi_c$ is $c$-optimal for the nonlinear heteroscedastic regression model (2.1) and (2.2) if and only if there exists a generalized inverse $G$ of the matrix $M(\xi_c, \theta)$ such that the inequality*

$$\frac{1}{\sigma^2(x,\theta)}\left(\frac{\partial\mu(x,\theta)}{\partial\theta}Gc\right)^2 + \frac{1}{2\sigma^4(x,\theta)}\left(\frac{\partial\sigma^2(x,\theta)}{\partial\theta}Gc\right)^2 \le c^T M^-(\xi_c,\theta)c \tag{3.2}$$

*holds for all $x \in \mathcal{X}$. Moreover, there is equality in (3.2) at any support point of the design $\xi_c$.*

In the following discussion, we will use the equivalence Theorem 3.1 to derive a geometric characterization of locally $c$-optimal designs in nonlinear regression models with Fisher information of the form (2.9), which generalizes the classical result of Elfving in an interesting direction. For this purpose, we define a generalized Elfving set by

$$\mathcal{R}_k = \mathrm{conv}\left\{\sum_{j=1}^k \varepsilon_j f_j(x,\theta) \,\middle|\, x \in \mathcal{X}; \sum_{j=1}^k \varepsilon_j^2 = 1\right\} \tag{3.3}$$

and obtain the following result.

THEOREM 3.3. *A design $\xi_c = \{x_r, p_r\}_{r=1}^m$ is locally $c$-optimal in a model with Fisher information matrix of the form (2.9) if and only if there exist constants $\gamma > 0$, $\varepsilon_{11}, \dots, \varepsilon_{1m}, \dots, \varepsilon_{k1}, \dots, \varepsilon_{km}$ satisfying*

$$\sum_{\ell=1}^k \varepsilon_{\ell r}^2 = 1, \qquad r = 1, \dots, m, \tag{3.4}$$

*such that the point $\gamma c \in \mathbb{R}^p$ lies on the boundary of the generalized Elfving set $\mathcal{R}_k$ defined in (3.3) and has the representation*

$$\gamma c = \sum_{r=1}^m p_r \left\{\sum_{\ell=1}^k \varepsilon_{\ell r} f_\ell(x_r,\theta)\right\} \in \partial\mathcal{R}_k. \tag{3.5}$$

PROOF. Assume that the design $\xi_c = \{x_r; p_r\}_{r=1}^m$ minimizes $c^T M^-(\xi,\theta)c$ or equivalently satisfies Theorem 3.1. We define $\gamma^2 = (c^T Gc)^{-1}, d = \gamma Gc \in \mathbb{R}^p$, then it follows from $c \in \mathrm{Range}(M(\xi,\theta))$ and the representation (2.9) that

$$\begin{aligned}\gamma c = M(\xi_c,\theta)d &= \sum_{r=1}^m p_r \left\{\sum_{\ell=1}^k f_\ell(x_r,\theta) f_\ell^T(x_r,\theta)\right\} d \\ &= \sum_{r=1}^m p_r \left\{\sum_{\ell=1}^k f_\ell(x_r,\theta)\varepsilon_{\ell r}\right\},\end{aligned}$$

where we have used the notation $\varepsilon_{\ell r} = f_\ell^T(x_r)d(\ell = 1,\dots,k; r = 1,\dots,m)$. By Theorem 3.1, there is equality in (3.1) for each support point $x_r$, which implies (observing the definition of $\gamma$ and $\varepsilon_{\ell r}$)

$$\sum_{\ell=1}^{k} \varepsilon_{\ell r}^2 = 1, \qquad r = 1,\dots,m.$$

Consequently, the conditions (3.4) and (3.5) are satisfied and it remains to show that $\gamma c \in \partial \mathcal{R}_k$. For this purpose, we note that we have by Cauchy's inequality and Theorem 3.1

$$\left(d^T\left\{\sum_{\ell=1}^{k} \varepsilon_\ell f_\ell(x,\theta)\right\}\right)^2 \le \sum_{\ell=1}^{k}(d^T f_\ell(x,\theta))^2 \sum_{\ell=1}^{k}\varepsilon_\ell^2 \le \sum_{\ell=1}^{k}\varepsilon_\ell^2 = 1 \tag{3.6}$$

for all $x \in \mathcal{X}$ and $\varepsilon_1,\dots,\varepsilon_\ell$ satisfying (3.4). Moreover, $\gamma c^T d = 1$ by the definition of the constant $\gamma$ and the vector $d$. Consequently, the vector $d$ defines a supporting hyperplane to the generalized Elfving set $\mathcal{R}_k$ at the point $\gamma c$, which implies $\gamma c \in \partial \mathcal{R}_k$

In order to prove the converse, assume that $\gamma c \in \partial \mathcal{R}_k$ and that (3.4) and (3.5) are satisfied. In this case, there exists a supporting hyperplane, say $d \in \mathbb{R}^p$, to the generalized Elfving set $\mathcal{R}_k$ at the boundary point $\gamma c$, that is,

$$\gamma c^T d = 1, \tag{3.7}$$

and the inequality

$$\left| d^T\left\{\sum_{\ell=1}^{k} \varepsilon_\ell f_\ell(x,\theta)\right\}\right| \le 1 \tag{3.8}$$

holds for all $x \in \mathcal{X}$, $\varepsilon_1,\dots,\varepsilon_k$ satisfying (3.4). Defining

$$\varepsilon_\ell(x) = \frac{d^T f_\ell(x,\theta)}{\sqrt{\sum_{j=1}^{k}(d^T f_j(x,\theta))^2}}, \qquad \ell = 1,\dots,k,$$

we have $\sum_{\ell=1}^{k} \varepsilon_\ell^2(x) = 1$ and obtain from the inequality (3.8)

$$\sum_{\ell=1}^{k}(d^T f_\ell(x))^2 = \left|\sum_{\ell=1}^{k} d^T\{\varepsilon_\ell(x) f_\ell(x,\theta)\}\right|^2 \le 1 \tag{3.9}$$

for all $x \in \mathcal{X}$. On the other hand, it follows from (3.5) and (3.7) that

$$1 = \gamma c^T d = \sum_{r=1}^{m} p_r d^T\left\{\sum_{\ell=1}^{k}\varepsilon_{\ell r} f_\ell(x_r,\theta)\right\},$$

which implies [using inequality (3.9)] for $x_1,\dots,x_m$,

$$d^T\left\{\sum_{\ell=1}^{k}\varepsilon_{\ell r} f_\ell(x_r,\theta)\right\} = 1, \qquad r = 1,\dots,m. \tag{3.10}$$

Consequently, we obtain using the Cauchy–Schwarz inequality

$$1 = \left(\sum_{\ell=1}^{k} \varepsilon_{\ell r} d^T f_\ell(x_r, \theta)\right)^2 \le \sum_{\ell=1}^{k} \varepsilon_{\ell r}^2 \sum_{\ell=1}^{k} (d^T f_\ell(x_r, \theta))^2 \le 1$$

for each $r = 1, \ldots, m$, which gives

$$\varepsilon_{\ell r} = \lambda_r d^T f_\ell(x_r, \theta), \qquad \ell = 1, \ldots, k, r = 1, \ldots, m, \tag{3.11}$$

for some constants $\lambda_1, \ldots, \lambda_m$. Now inserting (3.11) in (3.10) yields

$$1 = \lambda_r \sum_{\ell=1}^{k} (d^T f_\ell(x_r, \theta))^2 = \lambda_r, \qquad r = 1, \ldots, m,$$

and combining (3.11) and (3.5) gives

$$\gamma c = \sum_{r=1}^{m} p_r \left\{ \sum_{\ell=1}^{k} f_\ell(x_r, \theta) f_\ell^T(x_r, \theta) \right\} d = M(\xi_c, \theta) d.$$

It follows from Searle (1982) that there exists a generalized inverse of the matrix $M(\xi_c, \theta)$ such that $d = \gamma G c$ and we have from (3.9) for all $x \in \mathcal{X}$

$$\gamma^2 \sum_{\ell=1}^{k} (c^T G f_\ell(x))^2 \le 1. \tag{3.12}$$

By (3.7), we have $\gamma^2 = (c^T G c)^{-1} = (c^T M^{-1}(\xi_c, \theta) c)^{-1}$ and the inequality (3.12) reduces to (3.1). Consequently, the $c$-optimality of the design $\xi_c = \{x_r; p_r\}_{r=1}^m$ follows from Theorem 3.1. □

Note that in the case $k = 1$ the generalized Elfving set $\mathcal{R}_k$ defined in (3.3) reduces to the classical Elfving set (1.3) introduced by Elfving (1952). Similarly, Theorem 3.3 reduces to the classical Elfving theorem for models in which the Fisher information can be represented by $I(x, \theta) = f_1(x, \theta) f_1^T(x, \theta)$ [see Pukelsheim (1993), Chapter 2, or Studden (2005)]. In the following corollary, we specify the geometric characterization of Theorem 3.3 in the special case of nonlinear heteroscedastic regression models where we have $k = 2$ and the functions $f_1$ and $f_2$ are given by (2.10).

COROLLARY 3.4. *Consider the nonlinear heteroscedastic regression model defined by (2.1) and (2.2) with the additional assumption of a normal distribution, and define*

$$\mathcal{R} = \operatorname{conv}\left\{ \frac{\varepsilon_1}{\sigma(x, \theta)} \left(\frac{\partial \mu(x, \theta)}{\partial \theta}\right)^T + \frac{\varepsilon_2}{\sqrt{2}\sigma^2(x, \theta)} \left(\frac{\partial \sigma^2(x, \theta)}{\partial \theta}\right)^T \,\middle|\, x \in \mathcal{X}, \varepsilon_1^2 + \varepsilon_2^2 = 1 \right\}. \tag{3.13}$$

*A design $\xi_c = \{x_r; p_r\}_{r=1}^m$ is locally c-optimal if and only if there exist constants $\gamma > 0$ and $\varepsilon_{11}, \dots, \varepsilon_{1m}, \varepsilon_{21}, \dots, \varepsilon_{2m}$ satisfying $\varepsilon_{1r}^2 + \varepsilon_{2r}^2 = 1$ $(r = 1, \dots, m)$ such that the point*

$$\gamma c = \sum_{r=1}^m p_r \left\{ \frac{\varepsilon_{1r}}{\sigma(x_r,\theta)} \left( \frac{\partial \mu(x_r,\theta)}{\partial \theta} \right)^T + \frac{\varepsilon_{2r}}{\sqrt{2}\sigma^2(x_r,\theta)} \left( \frac{\partial \sigma^2(x_r,\theta)}{\partial \theta} \right)^T \right\}$$

*is a boundary point of the set $\mathcal{R}$.*

REMARK 3.5. Consider the case where mean and variance in the regression model are related by a known link function, say $\ell$, that is,

$$\sigma^2(x,\theta) = \ell(\mu(x,\theta)). \tag{3.14}$$

Under this assumption, a straightforward calculation shows that the Elfving set (3.13) reduces to

$$\begin{aligned} \mathcal{R} &= \operatorname{conv}\left\{ \frac{\varepsilon_1}{\sqrt{\ell(\mu(x,\theta))}} \left( \frac{\partial \mu(x,\theta)}{\partial \theta} \right)^T + \frac{\varepsilon_2 \ell'(\mu(x,\theta))}{\sqrt{2\ell^2(\mu(x,\theta))}} \left( \frac{\partial \mu(x,\theta)}{\partial \theta} \right)^T \,\Big|\, x \in \mathcal{X}, \varepsilon_1^2 + \varepsilon_2^2 = 1 \right\} \\ &= \operatorname{conv}\left\{ \left( \frac{\varepsilon_1}{\sqrt{\ell(\mu(x,\theta))}} + \frac{\varepsilon_2 \ell'(\mu(x,\theta))}{\sqrt{2\ell^2(\mu(x,\theta))}} \right) \left( \frac{\partial \mu(x,\theta)}{\partial \theta} \right)^T \,\Big|\, x \in \mathcal{X}, \varepsilon_1^2 + \varepsilon_2^2 = 1 \right\}. \end{aligned}$$

Now, using the values

$$\varepsilon_1 = \varepsilon \frac{1}{\sqrt{1+\omega^2(x)}}, \qquad \varepsilon_2 = \varepsilon \frac{\omega(x)}{\sqrt{1+\omega^2(x)}}, \qquad \omega(x) = \frac{\ell'(\mu(x,\theta))}{\sqrt{2\ell(\mu(x,\theta))}}$$

for some $\varepsilon \in \{-1, 1\}$ it can be easily shown that $\mathcal{R}_1 \subset \mathcal{R}$, where the set $\mathcal{R}_1$ is defined by

$$\mathcal{R}_1 = \operatorname{conv}\left\{ \varepsilon \sqrt{\frac{1}{\ell(\mu(x,\theta))} + \frac{1}{2}\left( \frac{\ell'(\mu(x,\theta))}{\ell(\mu(x,\theta))} \right)^2} \left( \frac{\partial \mu}{\partial \theta}(x,\theta) \right)^T \,\Big|\, x \in \mathcal{X}, \varepsilon \in \{-1, 1\} \right\},$$

and corresponds to a classical Elfving set of the form (1.3). Similarly, it follows observing the inequality

$$\left| \frac{\varepsilon_1}{\sqrt{\ell(\mu,(x,\theta))}} + \frac{\varepsilon_2 \ell'(\mu(x,\theta))}{\sqrt{2\ell^2(\mu(x,\theta))}} \right| \le \sqrt{\frac{1}{\ell(\mu(x,\theta))} + \frac{1}{2}\left( \frac{\ell'(\mu(x,\theta))}{\ell(\mu(x,\theta))} \right)^2}$$

whenever $\varepsilon_1^2 + \varepsilon_2^2 = 1$, that $\mathcal{R} \subset \mathcal{R}_1$, which implies $\mathcal{R} = \mathcal{R}_1$. Thus, for nonlinear heteroscedastic regression models, where the conditional expectation and variance are related by a known link function by (3.14), the Elfving set $\mathcal{R}$ reduces to an Elfving set of the form (1.3) and the classical result of Elfving (1952) can be used to characterize locally $c$-optimal designs.

**4. Examples.** In this section, we will illustrate the geometric characterization by two examples. In the first example, we discuss a nonlinear heteroscedastic regression model. Our second example considers a population model and indicates that an information matrix of the form (2.9) may also appear in other situations.

4.1. *A heteroscedastic Michaelis–Menten model.* As a first example, we consider a common model in enzyme kinetics with an heteroscedastic error structure, that is

$$Y_i = \frac{\theta_1 x_i}{\theta_2 + x_i} + \sqrt{e^{-\theta_3 x_i}}\varepsilon_i, \qquad i = 1, \ldots, n, \tag{4.1}$$

where errors are independent identically and normal distributed with mean 0 and variance $\sigma^2 > 0$. The parameters to be estimated are $\theta = (\theta_1, \theta_2, \theta_3)$. Clearly, $\mathrm{Var}(Y_i) = \tilde{\sigma}^2(x_i, \theta) := e^{-\theta_3 x_i}\sigma^2$, which means that the model is heteroscedastic. Some design issues for the Michaelis–Menten model have been discussed by López-Fidalgo and Wong (2002). In the following, we illustrate the geometric characterization of Section 3 by constructing optimal designs for estimating the minimum effective dose in the model (4.1).

For this purpose, we assume for the parameters $\theta = (3, 1.7, 0.1)$, $\sigma^2 = 1$. A straightforward calculation shows that

$$\frac{\partial \mu(x,\theta)}{\partial \theta} = \left(\frac{x}{\theta_2 + x}, \frac{-\theta_1 x}{(\theta_2 + x)^2}, 0\right)^T, \qquad \frac{\partial \tilde{\sigma}^2(x,\theta)}{\partial \theta} = (0, 0, -xe^{-\theta_3 x})^T,$$

which yields

$$f_1(x,\theta) = \frac{1}{\sigma e^{-\theta_3 x/2}}\left(\frac{x}{\theta_2 + x}, \frac{-\theta_1 x}{(\theta_2 + x)^2}, 0\right)^T, \tag{4.2}$$

$$f_2(x,\theta) = \frac{1}{\sqrt{2}\sigma^2 e^{-\theta_3 x}}(0, 0, -xe^{-\theta_3 x})^T. \tag{4.3}$$

The corresponding generalized Elfving space is as depicted in Figure 1. Suppose we are interested in estimating the minimum effective dose $x_{\min}$, that is, the smallest value of $x$ resulting in an expected value $E[Y|x] = E$. The solution of this equation is given by

$$x_{\min} = \frac{E\theta_2}{\theta_1 - E}.$$

An optimal design for estimating the minimum effective dose minimizes the variance of the estimate for $x_{\min}$. Consequently, if maximum likelihood is used to estimate $\theta$ and $x_{\min}$ is estimated by $\hat{x}_{\min} = E\hat{\theta}_2/(\hat{\theta}_1 - E)$ an optimal design for estimating the minimum effective dose is a locally $c$-optimal design problem for the vector

$$c = \left(\frac{-E\theta_2}{(\theta_1 - E)^2}, \frac{E}{\theta_1 - E}, 0\right)^T,$$

which is marked as the red line in Figure 1. Let $E = 1$ and the maximum possible observation $x_{\max} = 10$, that is $\mathcal{X} = [0, 10]$. From the figure, we obtain as optimal design for estimating the minimum effective dose the two point design

$$\xi_c = \begin{pmatrix} 1.1 & 10 \\ 0.967 & 0.033 \end{pmatrix}.$$

The two blue points denote the points of the form

$$\frac{\varepsilon_1}{\sigma(x,\theta)}\left(\frac{\partial\mu(x,\theta)}{\partial\theta}\right)^T + \frac{\varepsilon_2}{\sqrt{2}\sigma^2(x,\theta)}\left(\frac{\partial\sigma^2(x,\theta)}{\partial\theta}\right)^T \tag{4.4}$$

($x \in \mathcal{X}$, $\varepsilon_1^2 + \varepsilon_2^2 = 1$) which are used to represent the intersection of the Elfving set $\mathcal{R}$ defined by (3.13) with the line in direction of the vector $c$. The optimality of this design can be verified by Theorem 3.1.

4.2. *A random effect nonlinear regression model.* An information matrix of the form (2.9) appears also in the case of a nonlinear model with mixed effects and homoscedastic error structure, which has been intensively discussed in the toxicokinetics and pharmacokinetics literature [see, e.g., Beatty and Pigeorsch (1997) or Retout and Mentré (2003) among others]. To be precise, we consider a special case of a nonlinear mixed effects model which appears in population toxicokinetics, that is,

$$Y_i = f(x_i, b_i) + \varepsilon_i, \qquad i = 1, \ldots, n,$$

where the errors are independent identically and normally distributed with mean 0 and variance $\sigma^2 > 0$. The quantities $b_1, \ldots, b_n \sim \mathcal{N}(\theta, \Omega)$ denote here $p$-dimensional independent normally distributed random variables representing the effect of the corresponding subject under investigation [see Beatty and Pigeorsch (1997), Ette et al. (1995), Cayen and Black (1993)]. We assume that the random variables $b_1, \ldots, b_n$ and the vector $(\varepsilon_1, \ldots, \varepsilon_n)^T$ are independent. The variance of the random variable $Y_i$ can be approximated by

$$\operatorname{Var}(Y_i) \approx \tilde{\sigma}^2(x_i, \theta) := \frac{\partial f(x_i,\theta)}{\partial\theta}^T \Omega \frac{\partial f(x_i,\theta)}{\partial\theta} + \sigma^2.$$

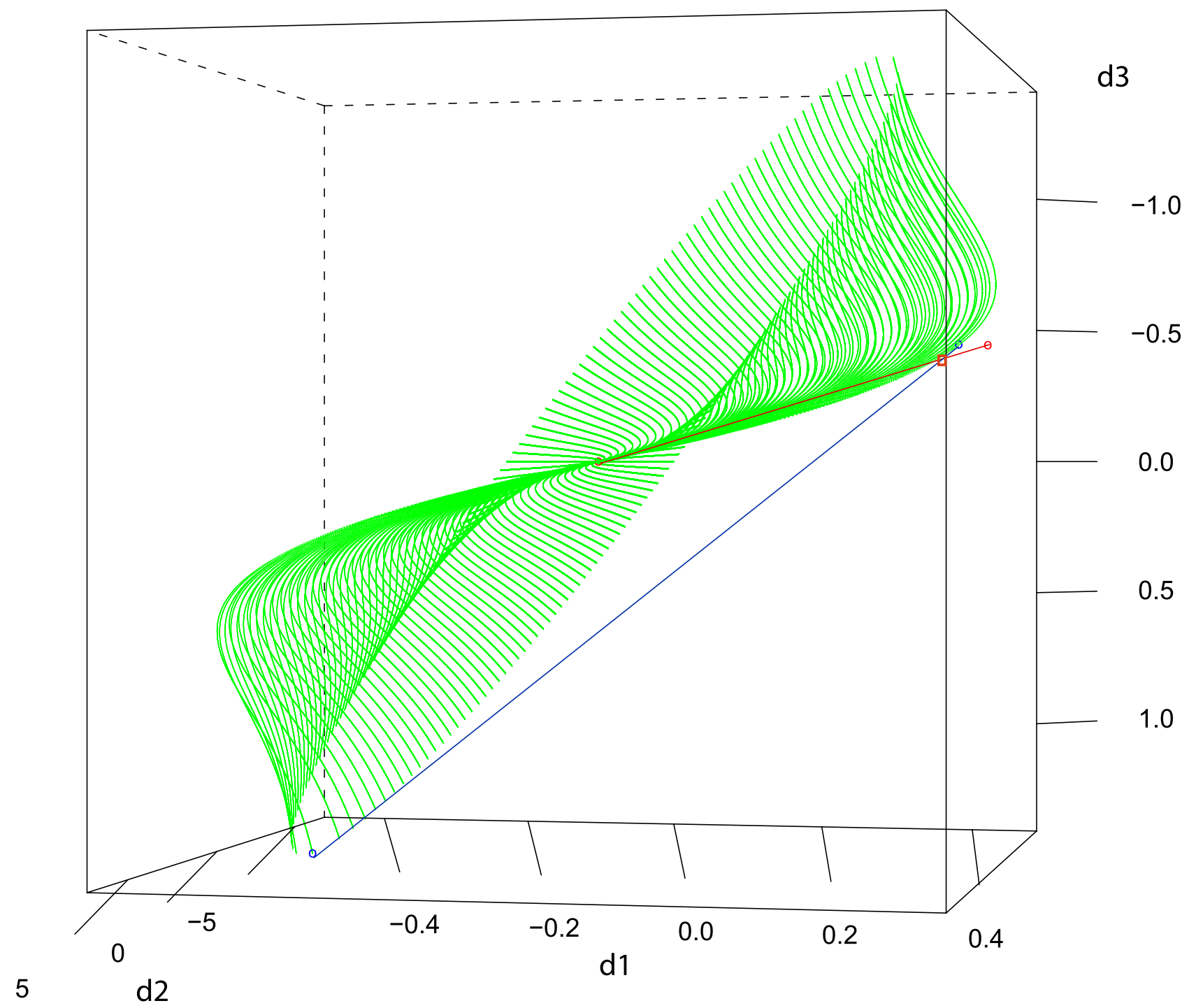

FIG. 1. *The Elfving space $\mathcal{R}$ defined in (3.3) for a Michaelis–Menten model (4.1) with heteroscedastic error structure. The functions $f_1$ and $f_2$ are given by (4.2) and (4.3), respectively. The blue line denotes a part of the boundary of the set $\mathcal{R}$.*

It now follows by similar arguments as in Retout and Mentré (2003) that the Fisher information matrix for the parameter $\theta$ at the point $x$ can be approximated by

$$I(x,\theta) = \frac{1}{\tilde{\sigma}^2(x,\theta)} \frac{\partial f(x,\theta)^T}{\partial\theta} \frac{\partial f(x,\theta)}{\partial\theta} + \frac{1}{2\tilde{\sigma}^4(x,\theta)} \frac{\partial \tilde{\sigma}^2(x,\theta)}{\partial\theta}^T \frac{\partial \tilde{\sigma}^2(x,\theta)}{\partial\theta},$$

which corresponds to the case (2.9) with $k = 2$ and

$$f_1(x,\theta) = \frac{1}{\tilde{\sigma}(x,\theta)} \frac{\partial f(x,\theta)}{\partial\theta}, \qquad f_2(x,\theta) = \frac{1}{\sqrt{2}\tilde{\sigma}^2(x,\theta)} \frac{\partial \tilde{\sigma}^2(x,\theta)}{\partial\theta}.$$

Consider, for example, the simple first-order elimination model

$$Y = b_1 e^{-b_2 x} + \varepsilon, \qquad x \in \mathcal{X} = [0,\infty), \tag{4.5}$$

which is widely used in pharmacokinetics [e.g., Rowland (1993)]. For the parameters, we assume

$$\theta = (30, 1.7), \qquad \Omega = \operatorname{diag}(\omega_1,\omega_2) = \operatorname{diag}(1, 0.1) \quad \text{and} \quad \sigma^2 = 0.04.$$

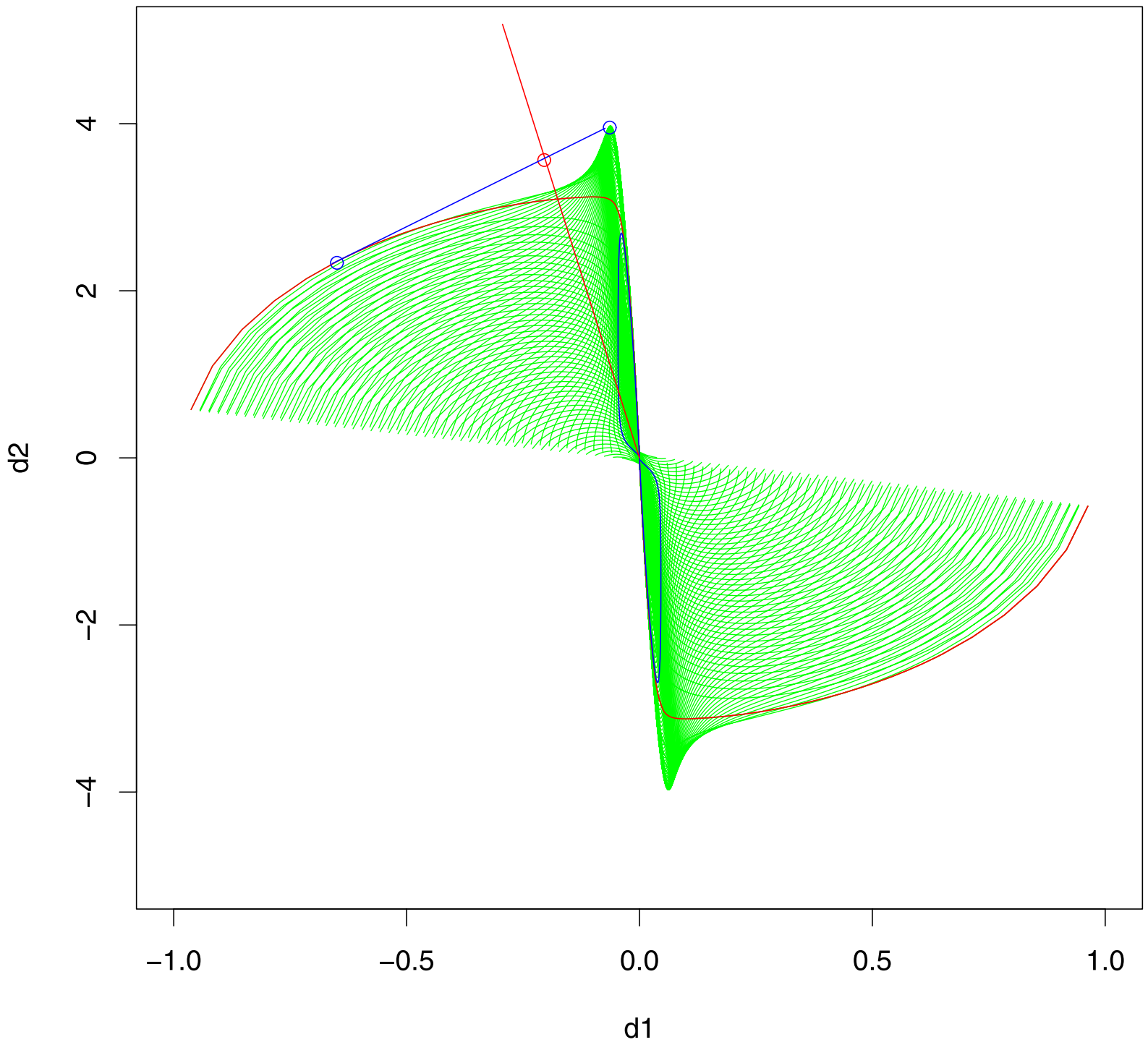

Fig. 2. *The Elfving space $\mathcal{R}$ defined in (3.3) for a simple population model (4.5). The functions $f_1$ and $f_2$ are given by (4.6) and (4.7), respectively. The blue line denotes a part of the boundary of the set $\mathcal{R}$.*

A straightforward calculation shows that

$$\frac{\partial f(x,\theta)}{\partial \theta} = (e^{-\theta_2 x}, -\theta_1 x e^{-\theta_2 x}),$$

$$\frac{\partial^2 f(x,\theta)}{\partial \theta_1 \, \partial \theta_2} = -x e^{-\theta_2 x}, \qquad \frac{\partial^2 f(x,\theta)}{\partial \theta_1^2} = 0, \qquad \frac{\partial^2 f(x,\theta)}{\partial \theta_2^2} = \theta_1 x^2 e^{-\theta_2 x},$$

which yields

$$f_1(x,\theta) = \frac{1}{\tilde{\sigma}(x,\theta)}(e^{-\theta_2 x}, -\theta_1 x e^{-\theta_2 x}), \tag{4.6}$$

$$f_2(x,\theta) = \frac{\sqrt{2}}{\tilde{\sigma}^2(x,\theta)} \sum_{m=1}^{2} \left( \frac{\partial^2 f(x,\theta)}{\partial \theta_m \, \partial \theta} \frac{\partial f(x,\theta)}{\partial \theta_m} \right) \omega_m. \tag{4.7}$$

The corresponding generalized Elfving space is as depicted in Figure 2. If we are interested in the optimal design for estimating the area under the curve, that is,

$$AUC = \int_0^\infty \theta_1 e^{-\theta_2 x}\, dx = \frac{\theta_1}{\theta_2},$$

we obtain a locally $c$-optimal design problem for the vector

$$c = (1/\theta_2, -\theta_1/\theta_2^2)^T,$$

which is marked as the red line in Figure 2. From this figure, we obtain as locally $c$-optimal design for the estimation of the area under the curve the two point design

$$\xi_c = \begin{pmatrix} 0.13 & 2.08 \\ 0.24 & 0.76 \end{pmatrix}.$$

The two blue points denote the points of the form (4.4) which are used to represent the intersection of the Elfving set $\mathcal{R}$ defined by (3.13) with the line in direction of the vector $c$. The optimality of this design can also be verified by Theorem 3.1.

**Acknowledgments.** The authors are grateful to Martina Stein who typed parts of this paper with considerable technical expertise. The authors are also grateful to two unknown referees for their constructive comments on an earlier version of this paper.

Fakultät für Mathematik
Ruhr-Universität Bochum
44780 Bochum
Germany
E-mail: holger.dette@rub.de

Medizinische Fakultät
Abteilung f. medizinische Informatik
Ruhr-Universität Bochum
44780 Bochum
Germany
E-mail: tim.holland-letz@rub.de